\documentclass{article}

\usepackage{arxiv}

\usepackage[utf8]{inputenc} % allow utf-8 input
\usepackage[T1]{fontenc}    % use 8-bit T1 fonts
\usepackage{hyperref}       % hyperlinks
\usepackage{url}            % simple URL typesetting
\usepackage{booktabs}       % professional-quality tables
\usepackage{amsfonts}       % blackboard math symbols
\usepackage{nicefrac}       % compact symbols for 1/2, etc.
\usepackage{microtype}      % microtypography
\usepackage{lipsum}
\usepackage{graphicx}
\graphicspath{ {./images/} }

\usepackage{authblk}
\usepackage{natbib}
\usepackage{subcaption}
\usepackage{tikz}
\usepackage{amsmath}
\usepackage[flushleft]{threeparttable}
\usepackage{multicol}
\usepackage{pgfplots}
\pgfplotsset{compat=1.18} 
\usepackage{bm}
\usepackage{algorithm}
\usepackage{algorithmic}
\newtheorem{theorem}{Theorem}[section]
\newtheorem{corollary}{Corollary}[theorem]
\newtheorem{lemma}[theorem]{Lemma}
\usepackage{enumerate}
\usepackage{float}

\title{A Modified Algorithm for Optimal Picker Routing \\
in a Single Block Warehouse}

\author[1]{George Dunn}
\author[2]{Hadi Charkhgard}
\author[3, 4]{Ali Eshragh}
\author[1]{Elizabeth Stojanovski}
\affil[1]{School of Information and Physical Sciences, University of Newcastle, NSW, Australia}
\affil[2]{Department of Industrial and Management Systems Engineering, University of South Florida, FL, USA}
\affil[3]{Carey Business School, Johns Hopkins University, MD, USA}
\affil[4]{International Computer Science Institute, University of California at Berkeley, CA, USA}

\date{}

\begin{document}
\maketitle
\begin{abstract}

The order picker routing problem involves finding the
optimal tour of a warehouse that collects all the required
items on a given pick list.
Ratliff and Rosenthal introduced a dynamic programming
algorithm for solving this problem in polynomial time
by sequentially adding edges inside and between
each aisle to construct a tour.
We provide a method where only transitions from one aisle
to the next are considered,
significantly reducing the number of stages in the
algorithm.

\end{abstract}

% ~~~~~~~~~~~~~~~~~~~~~~~~~~~~~~~~~~~~~~~~~~~~~~~~~~~~~~~~~~~~~~~~~~~~~~~~~~~~~~~~~
% Main Body
% ~~~~~~~~~~~~~~~~~~~~~~~~~~~~~~~~~~~~~~~~~~~~~~~~~~~~~~~~~~~~~~~~~~~~~~~~~~~~~~~~~

%%%%%%%%%%%%%%%%%%%%%%%%%%%%%%%%%%%%%%%%%%%%%%%%%%%%%%%%%%%%%%%%%%%%%%%%%%%%%%%
%%%%%%%%%%%%%%%%%%%%%%%%%%%%%%%%%%%%%%%%%%%%%%%%%%%%%%%%%%%%%%%%%%%%%%%%%%%%%%%
% INTRODUCTION
%%%%%%%%%%%%%%%%%%%%%%%%%%%%%%%%%%%%%%%%%%%%%%%%%%%%%%%%%%%%%%%%%%%%%%%%%%%%%%%
%%%%%%%%%%%%%%%%%%%%%%%%%%%%%%%%%%%%%%%%%%%%%%%%%%%%%%%%%%%%%%%%%%%%%%%%%%%%%%%

\section{Introduction}
\label{introduction}

Order picking is the process of collecting goods in a
warehouse to meet the requirements of the customer.
Given a pick list detailing the specific items needed,
these must be retrieved from their respective locations
within the warehouse and transported
to a depot for packing and shipping.
This has often been identified as a significant expense,
consuming up to 55\% of total operating costs,
highlighting the importance of optimizing related tasks
\citep{de2007design, tompkins2010facilities}. 
Determining the order in which the items should be collected
such that the distance traveled is minimized
is called the order picker routing problem.
This problem can be modeled as the well known
$\mathcal{NP}$-hard Travelling Salesman problem,
however \cite{ratliff1983order} presented a dynamic programming
algorithm for solving this problem in polynomial
time in the case of a single-block warehouse with a
parallel-aisle structure.

We present results showing that it is possible
to consider only the cross-aisle transitions
when constructing a tour,
removing a large number of steps
required to solve our updated algorithm.
For a warehouse with $m$ aisles,
the original method had $2m-1$ stages,
one for each aisle and each cross-aisle transition.
The algorithm we present solves the problem
with only $m$ stages.

Section \ref{ratliff} defines the picker routing problem
and provides background for the dynamic programming method of
\cite{ratliff1983order}.
We expand on these concepts
to produce a number of preliminary results in Section \ref{background}.
These allow us to present an updated algorithm
for solving the picker routing problem in a single block
warehouse in Section \ref{algorithm},
along with further modifications for a
rectangular layout.

%%%%%%%%%%%%%%%%%%%%%%%%%%%%%%%%%%%%%%%%%%%%%%%%%%%%%%%%%%%%%%%%%%%%%%%%%%%%%%%
%%%%%%%%%%%%%%%%%%%%%%%%%%%%%%%%%%%%%%%%%%%%%%%%%%%%%%%%%%%%%%%%%%%%%%%%%%%%%%%
% RELATED WORK
%%%%%%%%%%%%%%%%%%%%%%%%%%%%%%%%%%%%%%%%%%%%%%%%%%%%%%%%%%%%%%%%%%%%%%%%%%%%%%%
%%%%%%%%%%%%%%%%%%%%%%%%%%%%%%%%%%%%%%%%%%%%%%%%%%%%%%%%%%%%%%%%%%%%%%%%%%%%%%%

\section{Dynamic Programming Algorithm of Ratliff and Rosenthal}
\label{ratliff}

%%%%%%%%%%%%%%%%%%%%%%%%%%%%%%%%%%%%%%%%%%%%%%%%%%%%%%%%%%%%%%%%%%%%%%%%%%%%%%%
% Single Block Warehouse
%%%%%%%%%%%%%%%%%%%%%%%%%%%%%%%%%%%%%%%%%%%%%%%%%%%%%%%%%%%%%%%%%%%%%%%%%%%%%%%

\begin{figure}[ht]
\vskip 0.2in
\begin{center}
\centering
    \begin{subfigure}[b]{0.3\textwidth}
        \centering
        \resizebox{\linewidth}{!}{
            \begin{tikzpicture}[shorten >=1pt,draw=black!50]

    \draw[black, thin] (0, -1) -- (0 , 5.5) -- (7.5, 5.5) -- (7.5, -0.5)
    -- (1.5, -0.5) -- (1.5, -1) -- cycle;

    % Draw the input layer nodes
    \foreach \name / \x in {0,...,4}{
        \draw[black, thin] (1.5 * \x, 0) -- (1.5 * \x , 5);
        \draw[black, thin] (1.5 * \x + 0.5, 0) -- (1.5 * \x + 0.5 , 5);
        \draw[black, thin] (1.5 * \x + 1, 0) -- (1.5 * \x + 1 , 5);
        \draw[black, thin] (1.5 * \x + 1.5, 0) -- (1.5 * \x + 1.5, 5);
        \foreach \name / \y in {0,...,10}{
            \draw[black, thin] (1.5 * \x, 0.5 * \y) -- (1.5 * \x + 0.5, 0.5 * \y);
            \draw[black, thin] (1.5 * \x + 1, 0.5 * \y) -- (1.5 * \x + 1.5, 0.5 * \y);
            }
    }

    % \draw[black, thin] (0, -0.5) -- (0, -1) -- 
    % (1.5, -1) -- (1.5, -0.5) -- cycle;
    \node[align=center] at (0.75, -0.75) 
    {Depot};

    \node[align=center] at (3.75, -0.25) {Front Cross-Aisle};
    \node[align=center] at (3.75, 5.25) {Back Cross-Aisle};
    \node[align=center, rotate=90] at (3.75, 2.5) {Pick-Aisle};

    % \node[circle,fill=black,minimum size=5pt] (depot) at (3.75, -0.5) {};
    % \node[below of=depot, node distance=0.5cm] (d1) {Depot};

    \node[circle,fill=black,minimum size=3pt] (item1) at (0.25, 1.25) {};

    \node[circle,fill=black,minimum size=3pt] (item2) at (1.75, 4.25) {};

    \node[circle,fill=black,minimum size=3pt] (item3) at (2.75, 1.25) {};

    \node[circle,fill=black,minimum size=3pt] (item4) at (4.75, 3.75) {};

    \node[circle,fill=black,minimum size=3pt] (item5) at (5.75, 3.25) {};

    \node[circle,fill=black,minimum size=3pt] (item6) at (6.25, 0.75) {};

    \node[circle,fill=black,minimum size=3pt] (item7) at (7.25, 4.75) {};

\end{tikzpicture}
        }
        \caption{Warehouse}
        \label{fig:warehouse_fig}
    \end{subfigure}
    \begin{subfigure}[b]{0.3\textwidth}
    \centering
        \resizebox{\linewidth}{!}{
            \begin{tikzpicture}[shorten >=1pt,draw=black!50]

    % Lines

    \draw[black, thin] (0, -0.25) -- (6, -0.25);

    \draw[black, thin] (0, 5.25) -- (6, 5.25);

    \foreach \name / \y in {0,...,4}{
        \node[shape=circle,draw=black, minimum size=3pt, fill=white] (A-\name) at (1.5 * \y, -0.25) {};
        \node[shape=circle,draw=black, minimum size=3pt, fill=white] (B-\name) at (1.5 * \y , 5.25) {};

        \draw[black, thin] (A-\name) -- (B-\name);
        }

    % Items

    \node[circle, draw=black, minimum size=3pt, fill=black] (depot) at (0, -0.25) {};

    \node[circle, draw=black, minimum size=3pt, fill=black] (item1) at (0, 1.25) {};

    \node[circle,draw=black, minimum size=3pt, fill=black] (item2) at (1.5, 4.25) {};

    \node[circle,draw=black, minimum size=3pt, fill=black] (item3) at (1.5, 1.25) {};

    \node[circle,draw=black, minimum size=3pt, fill=black] (item4) at (4.5, 3.75) {};

    \node[circle,draw=black, minimum size=3pt, fill=black] (item5) at (4.5, 3.25) {};

    \node[circle,draw=black, minimum size=3pt, fill=black] (item6) at (6, 0.75) {};

    \node[circle,draw=black, minimum size=3pt, fill=black] (item7) at (6, 4.75) {};

    % Top and Bottom

    % \node[circle, draw=black, minimum size=3pt] (item1) at (0.75, -0.25) {};
    % \node[circle, draw=black, minimum size=3pt] (item1) at (0.75, 5.25) {};

\end{tikzpicture}
        }
        \caption{Graph Representation}   
        \label{fig:graph}
    \end{subfigure}
\caption{Standard Parallel Aisle Warehouse} 
\label{fig:warehouse}
\end{center}
\vskip -0.2in
\end{figure}
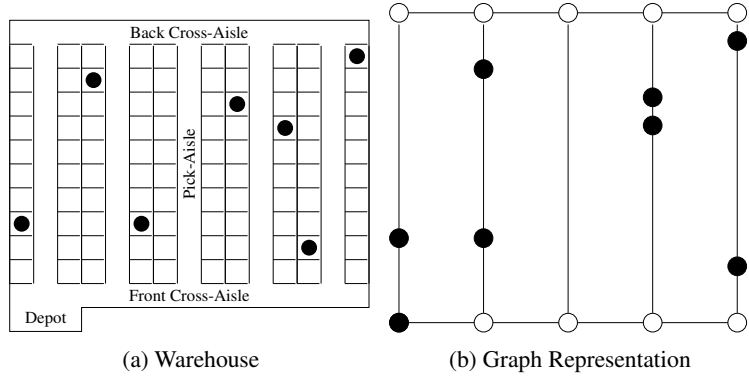

In this section, we provide an overview of the
picker routing problem and the corresponding
dynamic programming algorithm of \cite{ratliff1983order}.
We consider a single block warehouse with two horizontal cross-aisles and
$n \geq 1$ vertical pick aisles,
as seen in Figure \ref{fig:warehouse} where $n = 5$.
The aisles of the warehouse are assumed to be sufficiently
narrow such that the horizontal distance required to travel from one side
to the other is negligible.
The warehouse is represented by a graph $G$ with vertices $a_j$ and $b_j$
at the intersection of pick aisle $j=1,2,...,n$ and the
top and bottom cross-aisles, respectively.
A set of vertices $P = \{v_0, v_1,  ..., v_m\}$
represents the locations that must be visited,
with $v_0$ as the depot where the picker must start and finish a tour;
and $\{v_1, ..., v_m\}$ the products to be collected.
Only $v_0$ can be located at a $a_j$ or $b_j$ vertex,
all other vertices of $P$ are located within the picking aisles.

%%%%%%%%%%%%%%%%%%%%%%%%%%%%%%%%%%%%%%%%%%%%%%%%%%%%%%%%%%%%%%%%%%%%%%%%%%%%%%%
% Subgraph
%%%%%%%%%%%%%%%%%%%%%%%%%%%%%%%%%%%%%%%%%%%%%%%%%%%%%%%%%%%%%%%%%%%%%%%%%%%%%%%

%%% Tour Subgraph %%%
A subgraph $T \subset G$ is a tour subgraph if it contains
all vertices $v_i \in P$ and there exists an order picking
tour that uses each edge in $T$ exactly once.
The problem of finding an optimal order picking tour
can therefore be solved by finding a tour subgraph
with the minimum total edge length.

%%% Partial Tour Subgraph %%%
For $L \subset G$, a subgraph $T_j \subset L$ is a
Partial Tour Subgraph (PTS) if there exists
a subgraph $C_j \subset G - L$ such that
$T_j \cup C_j$ is a tour subgraph of $G$.
Two PTSs are said to be equivalent if any
completion of one is also a completion of the other.
It was shown that to show that two PTSs are
equivalence, it is sufficient to show that they
belong to the same equivalence class
defined by the degree parity of $a_j$ and $b_j$
and the number of connected components
in the graph.
There only exists seven possible
equivalence classes for a PTS:
\begin{equation*}
    \mathcal{E} = \{UU1C, 0E1C, E01C, EE1C, EE2C, 000C, 001C\}
\end{equation*}
where for each class,
the first two characters denote the
degree parity of $a_j$ and $b_j$ respectively
(zero 0, odd U or even E), and
the last two denote the connectivity
(0C, 1C or 2C).
The $000C$ class is possible only if none
of the aisles in $L_j$ contain any $v_i \in P$ vertices
and $001C$ is possible only if none of the
aisles in $G - L_j$ contain any.

%%% Construction of a Tour Subgraph %%%
\cite{ratliff1983order} showed that a tour subgraph
can be constructed by working
through all aisles from left to right,
first adding vertical edges within the aisle,
then horizontal edges between the current aisle and the
next.
The PTS that results from the addition of vertical edges in aisle
$j$ is called a $L_j^+$ PTS, and
the PTS that results from the addition of horizontal edges
from $j$ to $j+1$ is a $L_{j+1}^-$ PTS.

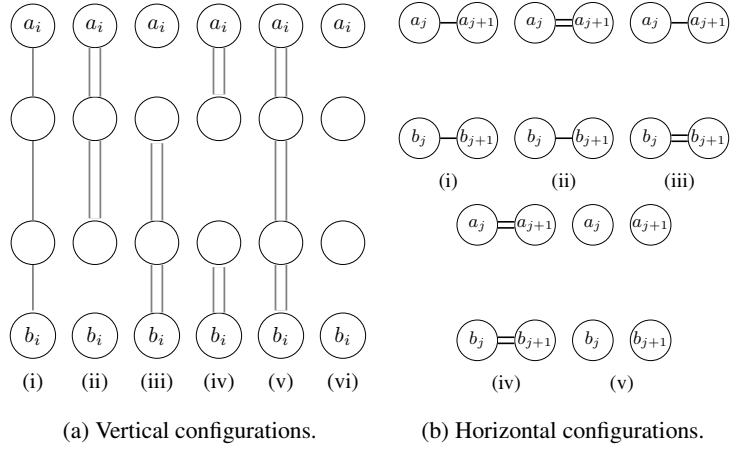
\begin{figure}[ht]
\vskip 0.1in
\begin{center}
\centering
    \begin{subfigure}[b]{0.30\textwidth}
        \centering
        \resizebox{\linewidth}{!}{
            % \begin{figure}
% \centering
%\vskip 0.2in
%\begin{center}
    
\begin{tikzpicture}[shorten >=1pt,->,draw=black!50, node distance=\layersep]

\tikzset{minimum size=20pt}

    % 1 Pass
    \node[shape=circle,draw=black] (a1) at (0, 5) {$a_i$};
    \node[shape=circle,draw=black] (d1) at (0, 3.5) { };
    \node[shape=circle,draw=black] (c1) at (0, 1.5) { };
    \node[shape=circle,draw=black] (b1) at (0, 0) {$b_i$}; 
    \draw[-, thick] (a1) -- (d1) -- (c1) -- (b1);

    \node (i) at (0, -0.75) {(i)};

    % Top
    \node[shape=circle,draw=black] (a2) at (1, 5) {$a_i$};
    \node[shape=circle,draw=black] (d2) at (1, 3.5) { };
    \node[shape=circle,draw=black] (c2) at (1, 1.5) { };
    \node[shape=circle,draw=black] (b2) at (1, 0) {$b_i$}; 
    \draw[-, thick, double, double distance between line centers=5pt] (a2) -- (d2) -- (c2);

    \node (ii) at (1, -0.75) {(ii)};

    % Bottom
    \node[shape=circle,draw=black] (a3) at (2, 5) {$a_i$};
    \node[shape=circle,draw=black] (d3) at (2, 3.5) { };
    \node[shape=circle,draw=black] (c3) at (2, 1.5) { };
    \node[shape=circle,draw=black] (b3) at (2, 0) {$b_i$}; 
    \draw[-, thick, double, double distance between line centers=5pt] (b3) -- (c3) -- (d3);

    \node (iii) at (2, -0.75) {(iii)};

    % Gap
    \node[shape=circle,draw=black] (a4) at (3, 5) {$a_i$};
    \node[shape=circle,draw=black] (d4) at (3, 3.5) { };
    \node[shape=circle,draw=black] (c4) at (3, 1.5) { };
    \node[shape=circle,draw=black] (b4) at (3, 0) {$b_i$}; 
    \draw[-, thick, double, double distance between line centers=5pt] (a4) -- (d4);
    \draw[-, thick, double, double distance between line centers=5pt] (b4) -- (c4);

    \node (iv) at (3, -0.75) {(iv)};

    % 2pass
    \node[shape=circle,draw=black] (a5) at (4, 5) {$a_i$};
    \node[shape=circle,draw=black] (d5) at (4, 3.5) { };
    \node[shape=circle,draw=black] (c5) at (4, 1.5) { };
    \node[shape=circle,draw=black] (b5) at (4, 0) {$b_i$}; 
    \draw[-, thick, double, double distance between line centers=5pt] (a5) -- (d5) -- (c5) -- (b5);

    \node (v) at (4, -0.75) {(v)};

    % None
    \node[shape=circle,draw=black] (a6) at (5, 5) {$a_i$};
    \node[shape=circle,draw=black] (d6) at (5, 3.5) { };
    \node[shape=circle,draw=black] (c6) at (5, 1.5) { };
    \node[shape=circle,draw=black] (b6) at (5, 0) {$b_i$}; 

    \node (vi) at (5, -0.75) {(vi)};

\end{tikzpicture}

%\vskip -0.1in

%\caption{Possible vertical edge configurations.}
%\label{fig:vertical_action}
%\end{center}
%\vskip -0.2in
%\end{figure}
        }
        \caption{Vertical configurations.}
        \label{fig:vertical_action}
    \end{subfigure}
    \begin{subfigure}[b]{0.29\textwidth}
    \centering
        \resizebox{\linewidth}{!}{
            %\begin{figure}
%\centering
%\vskip 0.2in
%\begin{center}

\def\v{1}
\def\h{0.5}
    
\begin{tikzpicture} % [shorten >=1pt,->,draw=black!50, node distance=\layersep]

\tikzset{minimum size=20pt}

    % 22
    \node[shape=circle,draw=black] (al1) at (2*\h, 2*\v) { };
    \node at (al1) {\small $a_j$}; 
    \node[shape=circle,draw=black] (ar1) at (4*\h, 2*\v) { };
    \node at (ar1) {\small $a_{j+1}$}; 
    \node[shape=circle,draw=black] (bl1) at (2*\h, 0*\v) { };
    \node at (bl1) {\small $b_j$}; 
    \node[shape=circle,draw=black] (br1) at (4*\h, 0*\v) { }; 
    \node at (br1) {\small $b_{j+1}$}; 
    \draw[-, thick] (al1) -- (ar1);
    \draw[-, thick] (bl1) -- (br1);

    \draw[-, thick, double, double distance between line centers=3pt] (al1) -- (ar1);
    \draw[-, thick, double, double distance between line centers=3pt] (bl1) -- (br1);

    \node (iv) at (3*\h, -0.75*\v) {(iv)};

    % 00
    \node[shape=circle,draw=black] (al2) at (6*\h, 2*\v) { };
    \node at (al2) {\small $a_j$}; 
    \node[shape=circle,draw=black] (ar2) at (8*\h, 2*\v) { };
    \node at (ar2) {\small $a_{j+1}$}; 
    \node[shape=circle,draw=black] (bl2) at (6*\h, 0*\v) { };
    \node at (bl2) {\small $b_j$}; 
    \node[shape=circle,draw=black] (br2) at (8*\h, 0*\v) { }; 
    \node at (br2) {\small $b_{j+1}$}; 

    \node (v) at (7*\h, -0.75*\v) {(v)};

    % 11
    \node[shape=circle,draw=black] (al3) at (0*\h, 5.5*\v) { };
    \node at (al3) {\small $a_j$}; 
    \node[shape=circle,draw=black] (ar3) at (2*\h, 5.5*\v) { };
    \node at (ar3) {\small $a_{j+1}$}; 
    \node[shape=circle,draw=black] (bl3) at (0*\h, 3.5*\v) { };
    \node at (bl3) {\small $b_j$}; 
    \node[shape=circle,draw=black] (br3) at (2*\h, 3.5*\v) { }; 
    \node at (br3) {\small $b_{j+1}$}; 
    \draw[-, thick] (al3) -- (ar3);
    \draw[-, thick] (bl3) -- (br3);

    \node (i) at (1*\h, 2.75*\v) {(i)};

    % 20
    \node[shape=circle,draw=black] (al4) at (4*\h, 5.5*\v) { };
    \node at (al4) {\small $a_j$}; 
    \node[shape=circle,draw=black] (ar4) at (6*\h, 5.5*\v) { };
    \node at (ar4) {\small $a_{j+1}$}; 
    \node[shape=circle,draw=black] (bl4) at (4*\h, 3.5*\v) { };
    \node at (bl4) {\small $b_j$}; 
    \node[shape=circle,draw=black] (br4) at (6*\h, 3.5*\v) { }; 
    \node at (br4) {\small $b_{j+1}$}; 
    \draw[-, thick] (al4) -- (ar4);
    \draw[-, thick] (bl4) -- (br4);

    \draw[-, thick, double, double distance between line centers=3pt] (al4) -- (ar4);

    \node (ii) at (5*\h, 2.75*\v) {(ii)};

    % 20
    \node[shape=circle,draw=black] (al5) at (8*\h, 5.5*\v) { };
    \node at (al5) {\small $a_j$}; 
    \node[shape=circle,draw=black] (ar5) at (10*\h, 5.5*\v) { };
    \node at (ar5) {\small $a_{j+1}$}; 
    \node[shape=circle,draw=black] (bl5) at (8*\h, 3.5*\v) { };
    \node at (bl5) {\small $b_j$}; 
    \node[shape=circle,draw=black] (br5) at (10*\h, 3.5*\v) { }; 
    \node at (br5) {\small $b_{j+1}$}; 
    \draw[-, thick] (al5) -- (ar5);
    \draw[-, thick] (bl5) -- (br5);

    \draw[-, thick, double, double distance between line centers=3pt] (bl5) -- (br5);

    \node (iii) at (9*\h, 2.75*\v) {(iii)};

\end{tikzpicture}

%\vskip -0.1in

%\caption{Possible horizontal edge configurations.}
%\label{fig:horizontal_action}
%\end{center}
%\vskip -0.2in
%\end{figure}
        }
        \caption{Horizontal configurations.}   
        \label{fig:horizontal_action}
    \end{subfigure}
\caption{Possible edge configurations in a optimal tour graph.} 
\label{fig:actions}
\end{center}
\vskip -0.2in
\end{figure}

%%% Actions %%%
For an optimum tour graph,
there are six possible vertical edge configurations
as shown in Figure \ref{fig:vertical_action}:
\begin{equation*}
    \mathcal{A} = \{1pass, top, bottom, gap, 2pass, none\}
\end{equation*}
a single traversal of the aisle ($1pass$),
entering and exiting from the back cross aisle ($top$),
entering and exiting from the bottom ($bottom$), 
traversing the aisle in a way that leaves the largest
section unconnected ($gap$),
traversing the aisle twice ($2pass$), and
not entering the aisle at all ($none$).
Note that $none$ is only valid if there are no
items to be collected within an aisle.
The resulting $L_j^+$ equivalence classes from the addition of
vertical edges to a $L_j^-$ PTS are shown in
Table \ref{tab:aisle_equivalence} in Appendix \ref{sec:rat_equ}.

For horizontal edges between aisles $j$ and $j+1$,
the five possible configurations are shown
in Figure \ref{fig:horizontal_action}:
\begin{equation*}
    \mathcal{H} = \{00, 20, 02, 22, 00\}
\end{equation*}
where the first and second numbers denote the number
of edges between the top and bottom cross-aisle vertices,
respectively.
The resulting $L_{j+1}^-$ equivalence classes from the addition of
vertical edges to a $L_j^+$ PTS are shown in
Table \ref{tab:cross_equivalence} in Appendix \ref{sec:rat_equ}.

The dynamic programming method works by
sequentially constructing PTSs,
adding all valid actions to each equivalence
class.
As we are interested in finding the optimal tour subgraph,
at each step, only the minimum PTS
for each equivalence class is stored.
When the minimum length $L_n^+$ PTS for each equivalence class
is found,
the optimal tour subgraph is the minimum of the
$E01C$, $0E1C$, $EE1C$ and $001C$ results.

%%%%%%%%%%%%%%%%%%%%%%%%%%%%%%%%%%%%%%%%%%%%%%%%%%%%%%%%%%%%%%%%%%%%%%%%%%%%%%%
%%%%%%%%%%%%%%%%%%%%%%%%%%%%%%%%%%%%%%%%%%%%%%%%%%%%%%%%%%%%%%%%%%%%%%%%%%%%%%%
% Background
%%%%%%%%%%%%%%%%%%%%%%%%%%%%%%%%%%%%%%%%%%%%%%%%%%%%%%%%%%%%%%%%%%%%%%%%%%%%%%%
%%%%%%%%%%%%%%%%%%%%%%%%%%%%%%%%%%%%%%%%%%%%%%%%%%%%%%%%%%%%%%%%%%%%%%%%%%%%%%%

\section{Preliminary Results}
\label{background}

%%%%%%%%%%%%%%%%%%%%%%%%%%%%%%%%%%%%%%%%%%%%%%%%%%%%%%%%%%%%%%%%%%%%%%%%%%%%%%%
% New Subgraph
%%%%%%%%%%%%%%%%%%%%%%%%%%%%%%%%%%%%%%%%%%%%%%%%%%%%%%%%%%%%%%%%%%%%%%%%%%%%%%%

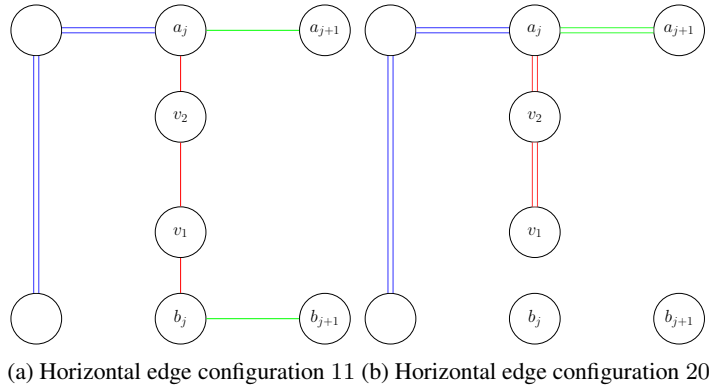
\begin{figure}[b]%[ht]
\vskip 0.1in
\begin{center}
\centering
    \begin{subfigure}[b]{0.28\textwidth}
        \centering
        \resizebox{\linewidth}{!}{
            \begin{tikzpicture} % [domain=0:4]
\node[circle,thick, draw=black, minimum size=50pt, font = {\LARGE}] (b1) at (0, 0) {$b_{j}$};
\node[circle,thick, draw=black, minimum size=50pt, font = {\LARGE}] (a1) at (0, 10) {$a_{j}$};
\node[circle,thick, draw=black, minimum size=50pt, font = {\LARGE}] (i1) at (0, 3) {$v_{1}$};
\node[circle,thick, draw=black, minimum size=50pt, font = {\LARGE}] (i2) at (0, 7) {$v_{2}$};

\node[circle,thick, draw=black, minimum size=50pt, font = {\LARGE}] (b2) at (5, 0) {$b_{j+1}$};
\node[circle,thick, draw=black, minimum size=50pt, font = {\LARGE}] (a2) at (5, 10) {$a_{j+1}$};

% Equivalence Class
\node[circle,thick, draw=black, minimum size=50pt] (a) at (-5, 10) {};
\node[circle,thick, draw=black, minimum size=50pt] (b) at (-5, 0) {};
% UU1C
% \draw[-, thick, draw=blue] (b1) -- (b) -- (a) -- (a1);
% E01C
\draw[-, thick, draw=blue, double, double distance between line centers=5pt] (a1) -- (a) -- (b);
% 0E1C
% \draw[-, thick, draw=blue, double, double distance between line centers=5pt] (b1) -- (b) -- (a);
% EE1C
% \draw[-, thick, draw=blue, double, double distance between line centers=5pt] (b1) -- (b) -- (a) -- (a1);
% EE2C
% \draw[-, thick, draw=blue, double, double distance between line centers=5pt] (b1) -- (b);
% \draw[-, thick, draw=blue, double, double distance between line centers=5pt] (a) -- (a1);

% Cross Actions
% 11
\draw[-, thick, draw=green] (a1) -- (a2);
\draw[-, thick, draw=green] (b1) -- (b2);
% 20
% \draw[-, thick, draw=green, double, double distance between line centers=5pt] (a1) -- (a2);
% 02
% \draw[-, thick, draw=green, double, double distance between line centers=5pt] (b1) -- (b2);
% 22
% \draw[-, thick, draw=green, double, double distance between line centers=5pt] (a1) -- (a2);
% \draw[-, thick, draw=green, double, double distance between line centers=5pt] (b1) -- (b2);

% Aisle Actions
% 1pass
\draw[-, thick, draw=red] (b1) -- (i1) -- (i2) -- (a1);
% Top
% \draw[-, thick, draw=red, double, double distance between line centers=5pt] (a1) -- (i2) -- (i1);
% Bottom
% \draw[-, thick, draw=red, double, double distance between line centers=5pt] (b1) -- (i1) -- (i2);
% Gap
% \draw[-, thick, draw=red, double, double distance between line centers=5pt] (b1) -- (i1);
% \draw[-, thick, draw=red, double, double distance between line centers=5pt] (a1) -- (i2);

\end{tikzpicture}
        }
        \caption{Horizontal edge configuration $11$}
        \label{fig:e01c_11}
    \end{subfigure}
    \begin{subfigure}[b]{0.28\textwidth}
    \centering
        \resizebox{\linewidth}{!}{
            \begin{tikzpicture} % [domain=0:4]
\node[circle,thick, draw=black, minimum size=50pt, font = {\LARGE}] (b1) at (0, 0) {$b_{j}$};
\node[circle,thick, draw=black, minimum size=50pt, font = {\LARGE}] (a1) at (0, 10) {$a_{j}$};
\node[circle,thick, draw=black, minimum size=50pt, font = {\LARGE}] (i1) at (0, 3) {$v_{1}$};
\node[circle,thick, draw=black, minimum size=50pt, font = {\LARGE}] (i2) at (0, 7) {$v_{2}$};

\node[circle,thick, draw=black, minimum size=50pt, font = {\LARGE}] (b2) at (5, 0) {$b_{j+1}$};
\node[circle,thick, draw=black, minimum size=50pt, font = {\LARGE}] (a2) at (5, 10) {$a_{j+1}$};

% Equivalence Class
\node[circle,thick, draw=black, minimum size=50pt] (a) at (-5, 10) {};
\node[circle,thick, draw=black, minimum size=50pt] (b) at (-5, 0) {};
% UU1C
% \draw[-, thick, draw=blue] (b1) -- (b) -- (a) -- (a1);
% E01C
\draw[-, thick, draw=blue, double, double distance between line centers=5pt] (a1) -- (a) -- (b);
% 0E1C
% \draw[-, thick, draw=blue, double, double distance between line centers=5pt] (b1) -- (b) -- (a);
% EE1C
% \draw[-, thick, draw=blue, double, double distance between line centers=5pt] (b1) -- (b) -- (a) -- (a1);
% EE2C
% \draw[-, thick, draw=blue, double, double distance between line centers=5pt] (b1) -- (b);
% \draw[-, thick, draw=blue, double, double distance between line centers=5pt] (a) -- (a1);

% Cross Actions
% 11
% \draw[-, thick, draw=green] (a1) -- (a2);
% \draw[-, thick, draw=green] (b1) -- (b2);
% 20
\draw[-, thick, draw=green, double, double distance between line centers=5pt] (a1) -- (a2);
% 02
% \draw[-, thick, draw=green, double, double distance between line centers=5pt] (b1) -- (b2);
% 22
% \draw[-, thick, draw=green, double, double distance between line centers=5pt] (a1) -- (a2);
% \draw[-, thick, draw=green, double, double distance between line centers=5pt] (b1) -- (b2);

% Aisle Actions
% 1pass
% \draw[-, thick, draw=red] (b1) -- (i1) -- (i2) -- (a1);
% Top
\draw[-, thick, draw=red, double, double distance between line centers=5pt] (a1) -- (i2) -- (i1);
% Bottom
% \draw[-, thick, draw=red, double, double distance between line centers=5pt] (b1) -- (i1) -- (i2);
% Gap
% \draw[-, thick, draw=red, double, double distance between line centers=5pt] (b1) -- (i1);
% \draw[-, thick, draw=red, double, double distance between line centers=5pt] (a1) -- (i2);

\end{tikzpicture}
        }
        \caption{Horizontal edge configuration $20$}   
        \label{fig:e01c_20}
    \end{subfigure}
\caption{Two $E01C$ $L_j^-$ PTSs (blue) with 
horizontal (green) and vertical edge (red) configurations.} 
\label{fig:equivalences}
\end{center}
\vskip 0.1in
\end{figure}

Our method exploits the fact that given a $L_j^-$ PTS after the
horizontal edge transition from aisle $j-1$ to aisle $j$,
regardless of the vertical edges added in aisle $j$,
there will only be the same five possible horizontal
configurations to transition from aisle $j$ to $j+1$.
We introduce a new subgraph,
$S_{j+1} \subset L_{j+1}^-$,
the result of adding a valid horizontal
edge configuration between aisles $j$ and $j+1$ to a $L_j^-$ PTS
(i.e. before the vertical edges have been added to aisle $j$).
We define the state of $S_{j+1}$ by
(the degree of $a_j$ and $b_j$, respectively,
the degree of $a_{j+1}$ and $b_{j+1}$, respectively,
and the number of connected components).
For the seven possible $L_j^-$ equivalence classes,
the $S_{j+1}$ states resulting from the five possible
horizontal configurations are presented in Table \ref{tab:new_equivalence}.
In Figure \ref{fig:equivalences} two examples are shown
for an arbitrary $E01C$ $L_j^-$ PTS represented by the blue edges.
Possible $S_{j+1}$ subgraphs are shown by the combination of
blue and green edges,
resulting in $(UU,UU,2C)$ and $(E0,E0,1C)$ states, respectively, and
the vertical edge configurations in red complete a
$L_{j+1}^-$ PTS.
The proofs of both lemmas presented in this section can be found in \cite{ratliff1983order}.
We use these to develop further conditions for producing a
minimal $L_{j+1}^-$ PTS from a $S_{j+1}$ subgraph.

\begin{table}[t]
\centering
%\begin{scriptsize}
\caption {$S_{j+1}$ subgraph states resulting from horizontal edge configurations}
\label{tab:new_equivalence}
\vskip 0.15in
\begin{center}
\begin{small}
\begin{sc}
\begin{threeparttable}
    \begin{tabular}{| c | c  c  c  c  c  |}
    \hline 
    $L_j^-$ & & \multicolumn{3}{c}{Horizontal Edge Configurations} & \\
    Class & $11$ & $20$ & $02$ & $22$ & $00$  \\
    \hline 
    $UU1C$ & $(EE,UU,1C)^{\ref{case:gap}}$ & $(UU,E0,1C)^{\ref{case:1pass}}$ & $(UU,0E,1C)^{\ref{case:1pass}}$ & $(UU,EE,1C)^{\ref{case:1pass}}$ & $(UU,00,1C)^{\ref{case:1pass}}$ \\
    $E01C$ & $(UU,UU,2C)^{\ref{case:1pass}}$ & $(E0,E0,1C)^{\ref{case:bottom}}$ & $(EE,0E,2C)^{\ref{case:2pass}}$ & $(EE,EE,2C)^{\ref{case:2gap}}$ & $(E0,00,1C)^{\ref{case:bottom}}$ \\
    $0E1C$ & $(UU,UU,2C)^{\ref{case:1pass}}$ & $(EE,E0,2C)^{\ref{case:2pass}}$ & $(0E,0E,1C)^{\ref{case:bottom}}$ & $(EE,EE,2C)^{\ref{case:2gap}}$ & $(0E,00,1C)^{\ref{case:bottom}}$ \\
    $EE1C$ & $(UU,UU,1C)^{\ref{case:1pass}}$ & $(EE,E0,1C)^{\ref{case:gap}}$ & $(EE,0E,1C)^{\ref{case:gap}}$ & $(EE,EE,1C)^{\ref{case:gap}}$ & $(EE,00,1C)^{\ref{case:gap}}$ \\
    $EE2C$ & $(UU,UU,2C)^{\ref{case:1pass}}$ & $(EE,E0,2C)^{\ref{case:2pass}}$ & $(EE,0E,2C)^{\ref{case:2pass}}$ & $(EE,EE,2C)^{\ref{case:2gap}}$ & $(EE,00,2C)^{\ref{case:2pass}}$ \\
    $000C$ & $(UU,UU,2C)^{\ref{case:1pass}}$ & $(E0,E0,1C)^{\ref{case:bottom}}$ & $(0E,0E,1C)^{\ref{case:bottom}}$ & $(EE,EE,2C)^{\ref{case:2gap}}$ & $(00,00,0C)^{\ref{case:depot}}$ \\
    $001C^a$ & $-$ & $-$ & $-$ & $-$ & $(00,00,1C)$ \\
    \hline
    \end{tabular}
    \begin{tablenotes}
        \small
        \item $^a$ This class is feasible only if there are no
        items to be picked to the right of aisle $j$.
    \end{tablenotes}
    %\vskip 0.2in
\end{threeparttable}
\end{sc}
\end{small}
\end{center}
%\end{scriptsize}
\end{table}

%%%%%%%%%%%%%%%%%%%%%%%%%%%%%%%%%%%%%%%%%%%%%%%%%%%%%%%%%%%%%%%%%%%%%%%%%%%%%%%
%%%%%%%%%%%%%%%%%%%%%%%%%%%%%%%%%%%%%%%%%%%%%%%%%%%%%%%%%%%%%%%%%%%%%%%%%%%%%%%
% Subgraph is a PTS Lemma
%%%%%%%%%%%%%%%%%%%%%%%%%%%%%%%%%%%%%%%%%%%%%%%%%%%%%%%%%%%%%%%%%%%%%%%%%%%%%%%
%%%%%%%%%%%%%%%%%%%%%%%%%%%%%%%%%%%%%%%%%%%%%%%%%%%%%%%%%%%%%%%%%%%%%%%%%%%%%%%

\subsection{Constructing a valid Partial Tour Subgraph}

The following lemma states the conditions required for a subgraph
to be a PTS.

\begin{lemma}
\label{lem:pts}
Necessary and sufficient conditions for $T_j \in L_j$ to be an $L_j$ PTS are
\begin{enumerate}[(i)]
    \item for all $v_i \in L_j$, the degree of $v_i$ is positive in $T_j$;
    \item every vertex in $T_j$, except possibly for $a_j$ and $b_j$,
    has an even degree; and 
    \item excluding vertices with zero degree,
    $T_j$ has either no connected component,
    a single connected component containing at least one of $a_j$ and $b_j$,
    or two connected components with $a_j$  in one component and $b_j$ in the other.
\end{enumerate}
\end{lemma}

As we start with a $L_j^-$ PTS,
we know these conditions are satisfied in all aisles $1,...j-1$.
Therefore, we only need to check that the vertical edges
added to aisle $j$ ensure the conditions are met in
aisles $j$ and $j+1$.

%%%%%%%%%%%%%%%%%%%%%%%%%%%%%%%%%%%%%%%%%%%%%%%%%%%%%%%%%%%%%%%%%%%%%%%%%%%%%%%
% All items are visited
%%%%%%%%%%%%%%%%%%%%%%%%%%%%%%%%%%%%%%%%%%%%%%%%%%%%%%%%%%%%%%%%%%%%%%%%%%%%%%%

%%%%%%%%%% No none in nonempty aisle %%%%%%%%%%

For the degree of all $v_i$ vertices in aisle $j$ to have a positive degree
(i.e. the items in that aisle are visited),
it is sufficient to ensure that any edge configuration other than $none$
is added to the aisle.

\begin{corollary}
    \label{cor:empty}
    The $none$ vertical edge configuration is not valid in
    an aisle containing $v_i$ vertices.
\end{corollary}

%%%%%%%%%%%%%%%%%%%%%%%%%%%%%%%%%%%%%%%%%%%%%%%%%%%%%%%%%%%%%%%%%%%%%%%%%%%%%%%
% Even degree
%%%%%%%%%%%%%%%%%%%%%%%%%%%%%%%%%%%%%%%%%%%%%%%%%%%%%%%%%%%%%%%%%%%%%%%%%%%%%%%

%%%%%%%%%% 1pass %%%%%%%%%%

Next, we need to ensure $a_j$ and $b_j$ have an even (including zero) degree.
From the possible vertical edge configurations in a minimum tour,
$1pass$ is the only one that adds an odd number of edges to
$a_j$ and $b_j$.
From this it is seen that if $a_j$ and $b_j$ have odd degrees,
the edge configuration $1pass$ is the only valid option.
If $a_j$ and $b_j$ have even degrees,
all configurations excluding $1pass$ meet this condition.

\begin{corollary}
\label{cor:1pass}
If the degrees of $a_j, b_j \in S_{j+1}$ are odd,
then $1pass$ is the only
valid vertical edge configuration in aisle $j$.
In all other states,
$1pass$ is not a valid configuration.
\end{corollary}

%%%%%%%%%%%%%%%%%%%%%%%%%%%%%%%%%%%%%%%%%%%%%%%%%%%%%%%%%%%%%%%%%%%%%%%%%%%%%%%
% Connected components
%%%%%%%%%%%%%%%%%%%%%%%%%%%%%%%%%%%%%%%%%%%%%%%%%%%%%%%%%%%%%%%%%%%%%%%%%%%%%%%

%%%%%%%%%% 2 Connected  %%%%%%%%%%

It is seen that if the graph has one connected component,
connecting the $v_i$ vertices in aisle $j$
and ensuring an even degree of $a_j$ and $b_j$ will
provide a $L_{j+1}^-$ PTS;
however, if $S_{j+1}$ has two connected components
without $a_{j+1}$ in one and
$b_{j+1}$ in the other,
then the vertical edge configuration in aisle $j$
must result in a one connected graph by connecting the two components.
%%%%%%%%%% 2pass or one 1pass  %%%%%%%%%%
As $L_j^-$ is a valid PTS,
if there are two connected components, then one is in $a_j$
and the other is in $b_j$.
Connecting $a_j$ with $b_j$ therefore connects the two components.
Of the possible edge configurations in aisle $j$,
$1pass$ and $2pass$ are the only two that connect
$a_j$ and $b_j$,
therefore connecting the two separate components to form one.

\begin{corollary}
\label{cor:pass}
If $S_{j+1}$ has two connected components without $a_{j+1}$ in one and
$b_{j+1}$ in the other,
$1pass$ and $2pass$ are the only valid vertical edge configurations
possible in aisle $j$.
\end{corollary}

%%%%%%%%%% Top and bottom  %%%%%%%%%%

If the parity of $a_j$ is zero after the horizontal edges added
to $S_{j+1}$,
then the parity of $a_{j+1}$ will also be zero.
In this scenario, the edge configurations $top$ and $gap$,
will create a new connected
component containing $a_j$ and no other cross-aisle vertices.
The result is therefore not a $L_{j+1}^-$ PTS,
as this component does not
contain $a_{j+1}$ or $b_{j+1}$.
Similarly, we see that $bottom$ and $gap$ produce a new connected component
if the degree of $b_j$ is zero.

\begin{corollary}
\label{cor:bottom}
If $a_j$ ($b_j$) has a degree of zero,
the $top$ ($bottom$) and $gap$ edge configurations are not valid.
\end{corollary}

%%%%%%%%%%%%%%%%%%%%%%%%%%%%%%%%%%%%%%%%%%%%%%%%%%%%%%%%%%%%%%%%%%%%%%%%%%%%%%%
%%%%%%%%%%%%%%%%%%%%%%%%%%%%%%%%%%%%%%%%%%%%%%%%%%%%%%%%%%%%%%%%%%%%%%%%%%%%%%%
% Equivalent PTS Lemma
%%%%%%%%%%%%%%%%%%%%%%%%%%%%%%%%%%%%%%%%%%%%%%%%%%%%%%%%%%%%%%%%%%%%%%%%%%%%%%%
%%%%%%%%%%%%%%%%%%%%%%%%%%%%%%%%%%%%%%%%%%%%%%%%%%%%%%%%%%%%%%%%%%%%%%%%%%%%%%%

\subsection{Finding the minimal Partial Tour Subgraph}

The next lemma is used to define the conditions for two PTSs
to be equivalent.

\begin{lemma}
\label{lem:eq}
Two $L_j$ PTSs $T_j^1$ and $T_j^2$ are equivalent if
\begin{enumerate}[(i)]
    \item $a_j$ has the same degree parity (even, odd, or zero)
    in both and $b_j$ has the same degree parity in both,
    and
    \item Excluding vertices with zero degree,
    both $T_j^1$ and $T_j^2$ have no connected component,
    both have a single connected component containing at least one of $a_j$ and $b_j$,
    or both have two connected components with $a_j$ in one component and $b_j$ in the other.
\end{enumerate}
\end{lemma}

%%%%%%%%%%%%%%%%%%%%%%%%%%%%%%%%%%%%%%%%%%%%%%%%%%%%%%%%%%%%%%%%%%%%%%%%%%%%%%%
% One connected
%%%%%%%%%%%%%%%%%%%%%%%%%%%%%%%%%%%%%%%%%%%%%%%%%%%%%%%%%%%%%%%%%%%%%%%%%%%%%%%

% Now that we have ensured a PTS...

As we are looking for an optimal tour,
we are only interested in the minimal PTS for each
equivalence class at each stage.
If multiple actions produce equivalent PTSs,
all but the minimal can be ignored.
For a $S_{j+1}$ subgraph,
the vertical edge configuration added to aisle $j$ will not
change the degrees of $a_{j+1}$ and $b_{j+1}$.
Therefore, the only way for the edge configuration to affect the
equivalence class is to change the number of connected components.

%%%%%%%%% Gap is the minimal out of Gap, Top and Bottom %%%%%%%%%%

If $S_{j+1}$ has one connected component and the degrees of $a_j$ and $b_j$ are even,
$top$, $bottom$, $gap$ and $2pass$ will all result in equivalent PTSs.
By definition, $gap$ always leaves the maximum distance within the aisle
unconnected and
therefore will always be equal to or less than the others. 

\begin{corollary}
\label{cor:gap}
If the degrees of $a_j$ and $b_j$ are even
and $S_{j+1}$ has one connected component,
then the $gap$ vertical edge configuration
in aisle $j$
will always be minimal.
\end{corollary}

%%%%%%%%% 2pass and gap %%%%%%%%%%

If $S_{j+1}$ has two connected components,
one containing $a_{j+1}$ and $b_{j+1}$ in the other
(this also implies that $a_j$ and $b_j$ are in
separate components),
and $a_j$ and $b_j$ have an even degree
(i.e. $1pass$ is not valid),
then $2pass$ will connect
the two parts and result in a different equivalence class
to $top$, $bottom$ and $gap$,
therefore must be considered as a possible minimal PTS
for a different class of equivalence.
As $top$, $bottom$ and $gap$ all produce the same equivalence
class,
we only need to consider the minimal,
which is always $gap$

\begin{corollary}
\label{cor:2gap}
    If $a_{j+1}, b_{j+1} \in S_{j+1}$ are contained in two separate components
    and $a_j, b_j \in S_{j+1}$ have even degree
    then $gap$ and $2pass$ are both possible minimal edge configurations
    for PTSs with different equivalence classes. 
\end{corollary}

%%%%%%%%%%%%%%%%%%%%%%%%%%%%%%%%%%%%%%%%%%%%%%%%%%%%%%%%%%%%%%%%%%%%%%%%%%%%%%%
%%%%%%%%%%%%%%%%%%%%%%%%%%%%%%%%%%%%%%%%%%%%%%%%%%%%%%%%%%%%%%%%%%%%%%%%%%%%%%%
% Empty Aisles
%%%%%%%%%%%%%%%%%%%%%%%%%%%%%%%%%%%%%%%%%%%%%%%%%%%%%%%%%%%%%%%%%%%%%%%%%%%%%%%
%%%%%%%%%%%%%%%%%%%%%%%%%%%%%%%%%%%%%%%%%%%%%%%%%%%%%%%%%%%%%%%%%%%%%%%%%%%%%%%

For a $S_{j+1}$ subgraph,
$none$, $top$, $bottom$ and $gap$ all result in the same
degree and connectivity.
In fact, by definition, the cost of each configuration is
the same if an aisle is empty; therefore, without loss of generality,
the $none$ configuration can be substituted for any of these actions.

If aisle $j$ is empty,
$2pass$ can still potentially lead to a graph with different connectivity,
therefore if $S_{j+1}$ has two connected components,
one containing $a_{j+1}$ and $b_{j+1}$ in the other
and $a_j$ and $b_j$ have an even degree ($1pass$ not valid),
then $2pass$ will connect
the two parts and result in a different class of equivalence
to $none$ and must be considered.

\begin{corollary}
    If aisle $j$ is empty and $S_{j+1}$ has zero or one connected components,
    then $none$ will always be optimal.
\end{corollary}

%%%%%%%%%%%%%%%%%%%%%%%%%%%%%%%%%%%%%%%%%%%%%%%%%%%%%%%%%%%%%%%%%%%%%%%%%%%%%%%
%%%%%%%%%%%%%%%%%%%%%%%%%%%%%%%%%%%%%%%%%%%%%%%%%%%%%%%%%%%%%%%%%%%%%%%%%%%%%%%
%%%% One Aisle %%%
%%%%%%%%%%%%%%%%%%%%%%%%%%%%%%%%%%%%%%%%%%%%%%%%%%%%%%%%%%%%%%%%%%%%%%%%%%%%%%%
%%%%%%%%%%%%%%%%%%%%%%%%%%%%%%%%%%%%%%%%%%%%%%%%%%%%%%%%%%%%%%%%%%%%%%%%%%%%%%%

Lastly, we consider the case of a warehouse with all $v_i \in P$ contained within
the aisle $j$.
The $L_j^-$ equivalence is $000C$ as there are no items to the left of aisle $j$.
As there are no items to the right of aisle $j$, the only valid
aisle $j$ to $j+1$ horizontal edge configuration is $00$.
The $top$, $bottom$ and $2pass$ configurations all
result in a valid $001C$ $L_{j+1}^-$ PTS.
The minimum tour graph will depend on the location of the depot.
If the depot is located at $b_j$ ($a_j$) then $top$ ($bottom$)
% involves completely traversing the aisle twice,
is the same as $2pass$, which will
never be minimal,
therefore $bottom$ ($top$) is optimal.

\begin{corollary}
\label{cor:depot}
If $S_{j+1}$ has no connected components and the degree of
$a_j$, $b_j$, $a_{j+1}$ and $b_{j+1}$ is zero,
then the minimal vertical arc configuration in aisle $j$
is $top$ if the depot is located at $a_j$ 
or $bottom$ if the depot is located at $b_j$.
\end{corollary}

%%%%%%%%%%%%%%%%%%%%%%%%%%%%%%%%%%%%%%%%%%%%%%%%%%%%%%%%%%%%%%%%%%%%%%%%%%%%%%%
%%%%%%%%%%%%%%%%%%%%%%%%%%%%%%%%%%%%%%%%%%%%%%%%%%%%%%%%%%%%%%%%%%%%%%%%%%%%%%%
% Algorithms
%%%%%%%%%%%%%%%%%%%%%%%%%%%%%%%%%%%%%%%%%%%%%%%%%%%%%%%%%%%%%%%%%%%%%%%%%%%%%%%
%%%%%%%%%%%%%%%%%%%%%%%%%%%%%%%%%%%%%%%%%%%%%%%%%%%%%%%%%%%%%%%%%%%%%%%%%%%%%%%

\section{Simplified Dynamic Programming Algorithm}
\label{algorithm}

%%%%%%%%%%%%%%%%%%%%%%%%%%%%%%%%%%%%%%%%%%%%%%%%%%%%%%%%%%%%%%%%%%%%%%%%%%%%%%%
%%%%%%%%%%%%%%%%%%%%%%%%%%%%%%%%%%%%%%%%%%%%%%%%%%%%%%%%%%%%%%%%%%%%%%%%%%%%%%%
% Cases
%%%%%%%%%%%%%%%%%%%%%%%%%%%%%%%%%%%%%%%%%%%%%%%%%%%%%%%%%%%%%%%%%%%%%%%%%%%%%%%
%%%%%%%%%%%%%%%%%%%%%%%%%%%%%%%%%%%%%%%%%%%%%%%%%%%%%%%%%%%%%%%%%%%%%%%%%%%%%%%

From the corollaries presented in Section \ref{background},
we develop six cases that cover
all possible $S_{j+1}$ states,
as labeled in\\
Table \ref{tab:new_equivalence}:
\begin{enumerate}[\bfseries \text{Case} 1]
    \item \label{case:1pass} $\mathbf{1pass}$: from Corollary \ref{cor:1pass}, 
    if the degrees of $a_j, b_j \in S_{j+1}^-$ are odd,
    then $1pass$ is the only valid vertical edge configuration.
    \item \label{case:2pass} $\mathbf{2pass}$: from Corollary \ref{cor:1pass} and \ref{cor:pass},
    if the degree of $a_j$ and $b_j$ are even and $S_{j+1}$ has two connected components
    without $a_{j+1}$ in one and
    $b_{j+1}$ in the other,
    then $2pass$ is the only valid vertical edge configuration.
    \item \label{case:gap} $\mathbf{gap}$: from Corollary \ref{cor:gap},
    if the degree of $a_j$ and $b_j$ are both even (but not zero)
    and $S_{j+1}$ has one connected component,
    then the $gap$ vertical edge configuration
    will always be the minimal.
    \item \label{case:bottom} $\mathbf{bottom}$ $\mathbf{(top)}$: from Corollary \ref{cor:bottom},
    if $a_j$ ($b_j$) has a degree of zero,
    then $bottom$ ($top$) is the only valid vertical edge configuration.
    \item \label{case:2gap} $\mathbf{gap}$ \textbf{or} $\mathbf{2pass}$: from Corollary \ref{cor:2gap},
    if $a_{j+1}, b_{j+1} \in S_{j+1}$ are contained in two separate components
    and $a_j, b_j \in S_{j+1}$ have even degree
    then $gap$ or $2pass$ are both possible minimal edge configurations
    but for different equivalence classes. 
    \item \label{case:depot} $\mathbf{bottom}$: from Corollary \ref{cor:depot},
    if $S_{j+1}$ has no connected components and the degree of
    $a_j$, $b_j$, $a_{j+1}$ and $b_{j+1}$ is zero,
    then the minimal vertical edge configuration in aisle $j$
    is $top$ if the depot is located at $a_j$ 
    or $bottom$ if the depot is located at $b_j$.
    Without loss of generality,
    we assume a depot located in the bottom cross
    aisle.
\end{enumerate}

From these cases we see that for all
possible $S_{j+1}$ states there are,
at most,
two vertical edge configurations for aisle $j$
that will produce a $L_{j+1}^-$ PTS.
If fact,
there are only four instances where two
configurations are possible,
and they are all in the $22$ horizontal
arc configurations.
By considering the $22$ configuration twice at each stage
(one time for each of the possible vertical configurations)
and all other configurations once,
we still only have six possible decisions per stage,
which is the same as the number of vertical decisions
in the original algorithm.
% Horizontal actions
For this we define a new horizontal configuration
set that includes $22^*$,
the combination of $2pass$ and $22$ configurations:
\begin{equation*}
    \mathcal{H}^* = \{11, 20, 02, 22, 22^*, 00\}
\end{equation*}
The resulting $L_{j+1}^-$ equivalence classes
(with the required aisle $j$ vertical configurations)
for applying all horizontal edge configurations
to all $L_j^-$ PTS equivalence classes is
shown in Table \ref{tab:2p_final_equivalence}.
%%%%%%%%%% Equivalences %%%%%%%%%%
\begin{table}[t]
\centering
%\begin{scriptsize}
\caption {$L_{j+1}^-$ Equivalence classes (and required vertical edges) for
each horizontal arc configuration}
\label{tab:2p_final_equivalence}
\vskip 0.15in
\begin{center}
\begin{small}
\begin{sc}
\begin{threeparttable}
    \begin{tabular}{| c | c  c  c  c  c  c  |}
    \hline 
    $L_j^-$ & & \multicolumn{4}{c}{Horizontal Arc Configurations} & \\
    Class & $11$ & $20$ & $02$ & $22$ & $22^*$ & $00$  \\
    \hline 
    $UU1C$ & $UU1C$ ($gap$) & $E01C$ ($1pass$) & $0E1C$ ($1pass$) & $EE1C$ ($1pass$) & $-$ & $001C$ ($1pass$) \\
    $E01C$ & $UU1C$ ($1pass$) & $E01C$ ($top$) & $0E1C$ ($2pass$) & $EE2C$ ($gap$) & $EE1C$ ($2pass$) & $001C$ ($top$) \\
    $0E1C$ & $UU1C$ ($1pass$) & $E01C$ ($2pass$) & $0E1C$ ($bottom$) & $EE2C$ ($gap$) & $EE1C$ ($2pass$) & $001C$ ($bottom$) \\
    $EE1C$ & $UU1C$ ($1pass$) & $E01C$ ($gap$) & $0E1C$ ($gap$) & $EE1C$ ($gap$) & $-$ & $001C$ ($gap$) \\
    $EE2C$ & $UU1C$ ($1pass$) & $E01C$ ($2pass$) & $0E1C$ ($2pass$) & $EE2C$ ($gap$) & $EE1C$ ($2pass$) & $001C$ ($2pass$) \\
    $000C$ & $UU1C$ ($1pass$) & $E01C$ ($top$) & $0E1C$ ($bottom$) & $EE2C$ ($gap$) & $EE1C$ ($2pass$) & $001C$ ($bottom$) \\
    $001C$ & $-$ & $-$ & $-$ & $-$ & $-$ & $000C$ ($none$) \\
    \hline
    \end{tabular}
    % \begin{tablenotes}
    %     \small
    %     \item $-$ Never a valid action.
    % \end{tablenotes}
    %\vskip 0.2in
\end{threeparttable}
\end{sc}
\end{small}
\end{center}
%\end{scriptsize}
\end{table}
%%%%%%%%%% Psuedocode %%%%%%%%%%
With this we can find an optimal
tour graph in a way similar to \cite{ratliff1983order},
but instead of constructing all $L_j^-$ and $L_j^+$ PTSs,
we instead only need to look at $L_j^-$.

To present the algorithm for finding a minimal
tour, 
we first define the set of equivalence classes
for valid PTSs in aisle $j$ as
$\mathcal{E}_j \subset \mathcal{E}$.
For a minimal $L_j^-$ PTS with equivalence $e$,
we let
$a_j^e \in \mathcal{A}$ be the previous vertical action,
$h_j^e \in \mathcal{H}^*$ be the previous horizontal action, and
$t_j^e \in \mathbb{R}_{\ge 0}$ is the total length of the minimal $L_j^-$ PTS.
% Costs
The costs of edge configurations are given by
$VerticalCost(a, j)$ for the cost of $a \in \mathcal{A}$ in aisle $j$ and
$HorizontalCost(h, j+1)$ for the cost of $h \in \mathcal{H}^*$ between
aisles $j$ and $j+1$.
Lastly,
$VerticleEdges(e, h)$ returns the required vertical
edge configuration when performing horizontal action $h \in \mathcal{H}^*$
on a $e \in \mathcal{E}$ PTS.
To find the vertical configuration in aisle $n$ that completes
the tour graph,
a dummy aisle $n+1$ is considered with no
$v_i$ vertices and a horizontal distance of
zero between $n$ and $n+1$.
The resulting algorithm is presented in Algorithm \ref{alg:dp},
where the optimal tour graph is given by the minimum $001C$ $L_{n+1}^-$ PTS.
The edge configurations of the graph can be
found by working backwards through the stages.

\begin{algorithm} %[t]
    \caption{Dynamic Programming Algorithm for Single-Block Rectangular Warehouse}
    \label{alg:dp}
\begin{algorithmic}
    % \STATE \textbf{Input:} number of Epochs E, 
    \STATE \textbf{input} warehouse $G$ with $n$ aisles
    % \STATE \textbf{create} dummy aisle $n+1$ with no items or horizontal distance
    \STATE \textbf{initialize} $\mathcal{E}_j \gets \{ \}$ \textbf{for} $j=1,2,...n+1$
    \STATE \textbf{let} $\mathcal{E}_1 \gets \{000C\}$, $t_1^{000C} = 0$
    \FOR{$j = 1,...,n$}
    
        \FORALL{$e_j \in \mathcal{E}_j$}
            \STATE $e \gets e_j$
            \STATE $t_{old} \gets t_j^e$
            \FORALL{$h \in H$}
                \STATE $a \gets VerticleEdges(e, h)$
                \STATE $t \gets t_{old} + VerticalCost(a, j) + HorizontalCost(h, j+1)$
                \STATE $e \gets NewEquivalence(e, h)$
                \IF{$e \notin \mathcal{E}_{j+1}$}
                    \STATE $\mathcal{E}_{j+1} \gets \mathcal{E}_{j+1} \cup \{e\}$
                    \STATE $a_{j+1}^e \gets a$
                    \STATE $h_{j+1}^e \gets h$
                    \STATE $t_{j+1}^e \gets t$
                \ELSIF{$t < t_{j+1}^e$}
                    \STATE $a_{j+1}^e \gets a$
                    \STATE $h_{j+1}^e \gets h$
                    \STATE $t_{j+1}^e \gets t$
                \ENDIF
            \ENDFOR
        \ENDFOR
    \ENDFOR
    % \STATE \textbf{return} $t_{n+1}^s | s = 001C\$

\end{algorithmic}
\end{algorithm}

%%%%%%%%%%%%%%%%%%%%%%%%%%%%%%%%%%%%%%%%%%%%%%%%%%%%%%%%%%%%%%%%%%%%%%%%%%%%%%%
% Rectangular Warehouse
%%%%%%%%%%%%%%%%%%%%%%%%%%%%%%%%%%%%%%%%%%%%%%%%%%%%%%%%%%%%%%%%%%%%%%%%%%%%%%%

\subsection{Rectangular Warehouse}

%%%%%%%%%% Warehouse Problem %%%%%%%%%%

This algorithm can be simplified even further when
considering the problem from the previous section
but with the added
assumption that $G$ has a rectangular shape, that is,
all aisles have the same length $d_{aisle}$ and
both cross-aisle segments between aisles $j$ and $j+1$ have the same
length $d_{j,j+1}^{cross}$ for $1 \leq j < n$.
For this warehouse shape,
\cite{revenant2024note} shows that there always
exists a minimum length tour that does not contain
a double cross of any aisle;
therefore, $2pass$ can be ignored.

%%%%%%%%%%%%%%%%%%%%%%%%%%%%%%%%%%%%%%%

By ignoring $2pass$,
all instances in Table \ref{tab:2p_final_equivalence}
that involve a $2pass$ configuration,
including the entire $22^*$ column,
as well as five other entries
that can be removed,
as they will never be optimal.
This means that there is now only one possible
vertical configuration for every combination
of horizontal edge and equivalence class,
a maximum of five valid decisions in each
stage of the algorithm.

%%%%%%%%%%%%%%%%%%%%%%%%%%%%%%%%%%%%%%%%%%%%%%%%%%%%%%%%%%%%%%%%%%%%%%%%%%%%%%%
% Complexity
%%%%%%%%%%%%%%%%%%%%%%%%%%%%%%%%%%%%%%%%%%%%%%%%%%%%%%%%%%%%%%%%%%%%%%%%%%%%%%%

\subsection{Complexity}

Both the original algorithm and our new method
require the calculation of cost coefficients,
which for a warehouse with $m$ locations to visit and $n$ aisles
was shown to run in $\mathcal{O} (m + n)$ time \citep{hessler2022note}.
However, the number of stages required in our algorithm is a significant
reduction,
now requiring only $n$ instead of the original $2n-1$,
while maintaining the same number of states per state.
We also see that for a rectangular warehouse,
our algorithm has a maximum of five valid decisions
per stage instead of six.

%%%%%%%%%%%%%%%%%%%%%%%%%%%%%%%%%%%%%%%%%%%%%%%%%%%%%%%%%%%%%%%%%%%%%%%%%%%%%%%
%%%%%%%%%%%%%%%%%%%%%%%%%%%%%%%%%%%%%%%%%%%%%%%%%%%%%%%%%%%%%%%%%%%%%%%%%%%%%%%
% CONCLUSION
%%%%%%%%%%%%%%%%%%%%%%%%%%%%%%%%%%%%%%%%%%%%%%%%%%%%%%%%%%%%%%%%%%%%%%%%%%%%%%%
%%%%%%%%%%%%%%%%%%%%%%%%%%%%%%%%%%%%%%%%%%%%%%%%%%%%%%%%%%%%%%%%%%%%%%%%%%%%%%%

\section{Conclusion}
\label{conclusion}

We have shown that when implementing the
algorithm of \cite{ratliff1983order} to construct a
minimal picker tour in a single block warehouse,
it is possible to skip the stages that decide the
aisle transitions.
Instead, the following cross-aisle transition will determine
with certainty what the previous configuration should be.
For a warehouse with $m$ aisle,
this results in an algorithm that has $m$ stages
as opposed to $2m-1$,
with the number of states per stage remaining the same.
Modifications are also presented for problems in a rectangular warehouse
that further reduce the number of decisions available per state.

As the original algorithm provided the foundation
for extensions to two-block \citep{roodbergen2001routing}
and multi-block warehouses \citep{scholz2016new, pansart2018exact},
the results of this paper have the potential to improve
these variations in a similar way.
Applications to unconventional warehouse configurations such as
fish-bone \citep{ccelk2014order}
and chevron \citep{masae2020optimal} also
provide avenues for further research.

% ~~~~~~~~~~~~~~~~~~~~~~~~~~~~~~~~~~~~~~~~~~~~~~~~~~~~~~~~~~~~~~~~~~~~~~~~~~~~~~~~~
% Acknowledgements
% ~~~~~~~~~~~~~~~~~~~~~~~~~~~~~~~~~~~~~~~~~~~~~~~~~~~~~~~~~~~~~~~~~~~~~~~~~~~~~~~~~

\section*{Acknowledgements}

George was supported by an Australian Government Research Training
Program (RTP) Scholarship.

% ~~~~~~~~~~~~~~~~~~~~~~~~~~~~~~~~~~~~~~~~~~~~~~~~~~~~~~~~~~~~~~~~~~~~~~~~~~~~~~~~~
% Bibliography
% ~~~~~~~~~~~~~~~~~~~~~~~~~~~~~~~~~~~~~~~~~~~~~~~~~~~~~~~~~~~~~~~~~~~~~~~~~~~~~~~~~

\bibliographystyle{apalike}

% ~~~~~~~~~~~~~~~~~~~~~~~~~~~~~~~~~~~~~~~~~~~~~~~~~~~~~~~~~~~~~~~~~~~~~~~~~~~~~~~~~
% Appendix
% ~~~~~~~~~~~~~~~~~~~~~~~~~~~~~~~~~~~~~~~~~~~~~~~~~~~~~~~~~~~~~~~~~~~~~~~~~~~~~~~~~

%%%%%%%%%%%%%%%%%%%%%%%%%%%%%%%%%%%%%%%%%%%%%%%%%%%%%%%%%%%%%%%%%%%%%%%%%%%%%%%
%%%%%%%%%%%%%%%%%%%%%%%%%%%%%%%%%%%%%%%%%%%%%%%%%%%%%%%%%%%%%%%%%%%%%%%%%%%%%%%
% APPENDIX
%%%%%%%%%%%%%%%%%%%%%%%%%%%%%%%%%%%%%%%%%%%%%%%%%%%%%%%%%%%%%%%%%%%%%%%%%%%%%%%
%%%%%%%%%%%%%%%%%%%%%%%%%%%%%%%%%%%%%%%%%%%%%%%%%%%%%%%%%%%%%%%%%%%%%%%%%%%%%%%
\newpage
\appendix
% \onecolumn
\section{Appendix}

%%%%%%%%%%%%%%%%%%%%%%%%%%%%%%%%%%%%%%%%%%%%%%%%%%%%%%%%%%%%%%%%%%%%%%

\subsection{Resulting equivalent classes when adding
valid edge configurations}
\label{sec:rat_equ}

\begin{table}[H]
\centering
\caption {Resulting $L_j^+$ Equivalence Classes from Adding
Vertical Arc Configurations \citep{ratliff1983order}.}
\label{tab:aisle_equivalence}
\vskip 0.15in
\begin{center}
\begin{small}
\begin{sc}
\begin{threeparttable}
%\begin{scriptsize}

    \begin{tabular}{| c | c  c  c  c  c  c  |}
    \hline 
    $L_j^-$ & & \multicolumn{4}{c}{Vertical Edge Configurations} & \\
    % \hline 
    Class & $1pass$ & $top$ & $bottom$ & $gap$ & $2pass$ & $none$$^a$ \\
    \hline 
    $UU1C$ & $EE1C$ & $UU1C$ & $UU1C$ & $UU1C$ & $UU1C$ & $UU1C$ \\
    $E01C$ & $UU1C$ & $E01C$ & $EE2C$ & $EE2C$ & $EE1C$ & $E01C$ \\
    $0E1C$ & $UU1C$ & $EE2C$ & $0E1C$ & $EE2C$ & $EE1C$ & $0E1C$ \\
    $EE1C$ & $UU1C$ & $EE1C$ & $EE1C$ & $EE1C$ & $EE1C$ & $EE1C$ \\
    $EE2C$ & $UU1C$ & $EE2C$ & $EE2C$ & $EE2C$ & $EE1C$ & $EE1C$ \\
    $000C$$^b$ & $UU1C$ & $E01C$ & $0E1C$ & $EE2C$ & $EE1C$ & $000C$ \\
    $001C$$^c$ & $-$ & $-$ & $-$ & $-$ & $-$ & $001C$ \\
    \hline
    \end{tabular}
    \begin{tablenotes}
        \small
        \item $^a$ this is not a feasible configuration if
        there is any item to be picked to the left of aisle $j$.
        \item $^b$ This class can occur only if there are no
        items to be picked to the left of aisle $j$.
        \item $^c$ This class is feasible only if there are no
        items to be picked to the right of aisle $j$.
        \item $^d$ Could never be optimal.
    \end{tablenotes}
    \vskip -0.1in
\end{threeparttable}
\end{sc}
\end{small}
\end{center}
%\end{scriptsize}
\end{table}

\begin{table}[H]
\centering
%\begin{scriptsize}
\caption {Resulting $L_{j+1}^-$ Equivalence Classes from Adding
Horizontal Edge Configurations \citep{ratliff1983order}.}
\label{tab:cross_equivalence}
\vskip 0.15in
\begin{center}
\begin{small}
\begin{sc}
\begin{threeparttable}
    \begin{tabular}{| c | c  c  c  c  c  |}
    \hline 
    $L_j^+$ & \multicolumn{5}{c |}{Horizontal Edge Configurations} \\
    Class & $11$ & $20$ & $02$ & $22$ & $00$ \\
    \hline 
    $UU1C$ & $UU1C$ & $-^a$ & $-^a$ & $-^a$ & $-^a$ \\
    $E01C$ & $-^a$ & $E01C$ & $-^b$ & $EE2C$ & $001C$ \\
    $0E1C$ & $-^a$ & $-^b$ & $0E1C$ & $EE2C$ & $001C$ \\
    $EE1C$ & $-^a$ & $E01C$ & $0E1C$ & $EE1C$ & $001C$ \\
    $EE2C$ & $-^a$ & $-^b$ & $-^b$ & $EE2C$ & $-^b$ \\
    $000C$ & $-^c$ & $-^c$ & $-^c$ &$-^c$ & $001C$ \\
    $001C$ & $-^b$ & $-^b$ & $-^b$ & $-^b$ & $001C$ \\
    \hline
    \end{tabular}
    \begin{tablenotes}
        \small
        \item $^a$ The degrees of $a_i$ and $b_i$ are odd.
        \item $^b$ No completion can connect the graph.
        \item $^c$ Would never be optimal.
    \end{tablenotes}
    \vskip -0.2in
\end{threeparttable}
\end{sc}
\end{small}
\end{center}
%\end{scriptsize}
\end{table}

% ~~~~~~~~~~~~~~~~~~~~~~~~~~~~~~~~~~~~~~~~~~~~~~~~~~~~~~~~~~~~~~~~~~~~~~~~~~~~~~~~~

\end{document}